\documentstyle[12pt]{article}
\newtheorem{theorem}{Theorem}

\newtheorem{lemma}[theorem]{Lemma}

\newtheorem{definition}[theorem]{Definition}

\newtheorem{remark}{Remark}

\font\bb=msbm10
\begin{document}

\author{B. Shklyar}
\title{{\bf On Attainable Set and Controllability for Abstract Control
Problem with
Unbounded Input Operator}}
\maketitle

\begin{abstract}
For linear evolution control system described by
$\dot{x}=Ax(t)+Bu(t),x(0)=x_{0}$ 
($A$ generates a strongly continuous semigroup $\{S(t)\}_{t\ge 0}$ on a
Banach space 
$X$; $B$ is a linear unbounded  operator), the attainable set $K\left(
t\right)$ set 
is studied. Conditions of the independence of $t$ for its closure
$\overline{K(t)}$ 
are established. Controllability conditions for some classes of evolution
systems are
obtained.
\end{abstract}

{\bf Keywords.} Attainable sets, controllability, abstract evolution
equations, linear hereditary systems.

{\bf AMS MOS subject classification.} 35R30, 35J20, 65M30.

\section{Introduction.}

We consider a system described by linear abstract differential equation of
evolution type 
\begin{eqnarray}
\dot{x}(t) &=&Ax(t)+Bu(t)  \label{e1} \\
x(0) &=&x_{0}  \label{e2}
\end{eqnarray}
where $X,U$ be Banach spaces, $x(t)\in X$ is the current state, $x_{0}\in X$
is the initial state, $u(t)\in U,u(\cdot )\in L_{2}([0,t_{1}],U)$ is the
control, $A$ is a linear operator generating a strongly continuous semigroup 
$\{S(t)\}_{t\ge 0}$ of operators in the class $C_{0}$ $;$ $B:U\rightarrow X$
is a linear possibly unbounded operator.

Let $x\left( t,x_{0},u\left( \cdot \right) \right) $ be the weak solution of
equation (\ref{e1})-(\ref{e2}), corresponding to the control $u\left( \cdot
\right) $.

\begin{definition}
A state $x\in X$ is said to be attainable in the time $t$ from the origin,
if there exists an admissible control $u\left( \tau \right) ,0\leq \tau \leq
t$, such that $x\left( t,0,u\left( \cdot \right) \right) =x$.
\end{definition}

\begin{definition}
The set $K\left( t\right) \,$of all states $x\in X$ attainable in the time $t
$ from the origin is said to be attainable set for equation (1).
\end{definition}

Thus the {\it attainable set} $K(t)$ for equation (1) is defined as 
\begin{equation}  \label{e3}
K(t)=\left\{ x\in X:\exists u(\cdot )\in L_2\left( [0,t],U\right) ,x=x\left(
t,x_0,u\left( \cdot \right) \right) .\,\right\} .
\end{equation}

It is our main purpose here to study the properties of the set $K\left(
t\right) ,$ and especially, to establish conditions for its independence of $%
t$ for sufficiently large $t.$

These conditions has been established in author$^{\prime }$s paper \cite[ ]
{Shk95} for the equation (\ref{e1}) with bounded operator $B,$ so this work
can be considered as the continuation of \cite[ ]{Shk95}. The reason to
expand the results of \cite[ ]{Shk95} to equation (\ref{e1}) with unbounded
operator $B$ is the existence of the large classes of infinite dimensional
control systems adequately described by equation (\ref{e1}) with unbounded
operator $B$. Some of these classes are:

\begin{itemize}
\item  partial differential equations with boundary control;

\item  functional differential equations with delays both in state and in
control variables.
\end{itemize}

The importance of equation (\ref{e1}) with unbounded operator $B$ (both from
theoretical and from practical point of view) has been recognized by many
authors. We will use the functional analytic approach developed by Salamon 
\cite{Slmn87}.


\section{Preliminaries.}

Denote by $\sigma $ the spectrum of the operator $A$. Let $\mu \notin \sigma 
$. We will consider the spaces $W$ and $V$ defined as follows \cite{DPrt84},%
\cite{Ngl86},\thinspace \cite{Weiss89}:

$W$ is the domain $D\left( A\right) $ of the operator $A$ with the norm $%
\left\| x\right\| _\mu =\left\| \left( \mu I-A\right) x\right\| $;

$V$ is the closure of $X$ with respect to the norm $\left\| x\right\| _{-\mu
}=\left\| \left( \mu I-A\right) ^{-1}x\right\| .$

Obviously $W\subset X\subset V$.

It is known that for each $\mu _{1},\mu _{2},$ $\mu _{1}\neq \mu _{2}$ the
norm $\left\| \cdot \right\| _{\mu _{1}}$ is equivalent to the norm $\left\|
\cdot \right\| _{\mu _{2}}$, the norm $\left\| \cdot \right\| _{-\mu _{1}}$
is equivalent to the norm $\left\| \cdot \right\| _{-\mu _{2}}$, and all $%
\left\| \cdot \right\| _{\mu }$ are equivalent to the appropriate graph norm
on $D\left( A\right) $, so the Banach spaces $W$ and $V$ do not depend of $%
\mu $ \cite{Weiss89}.

We assume

\begin{enumerate}
\item  the unbounded operator $B$ is bounded as operator from $U$ to $V$ ;

\item  the operator $\Phi \left( t\right) :L_{2}\left( \left[ 0,t\right]
,U\right) \rightarrow X$ defined by the formula $\Phi \left( t\right)
u\left( \cdot \right) =\int_{0}^{t}S\left( t-\tau \right) Bu\left( \tau
\right) d\tau $ is bounded for each $t\geq 0$\footnote{%
It is easy to show, that $\Phi \left( t\right) u\left( \cdot \right)
=x\left( t,0,u\left( \cdot \right) \right) $.}.
\end{enumerate}

\medskip If $x\in X$ and $f\in X^{*}$, we will write $(x,f)$ instead of $%
f(x) $. The upper superscript ${\rm {}}^{{\rm T}}$ denotes transposition.

As usual ${\hbox{\bb R}}$ is the set of real and ${\hbox{\bb C}}$ the set of
complex numbers.

\medskip For any set $K\subset X$ we denote by $\overline{K}$ the closure of 
$K$ with respect to the uniform topology of $X$ and by $K^{\perp }$ the set $%
\{y\in X^{*}:\ (x,y)=0\ \forall x\in K\}$.

\medskip We assume the operator $A$ to have the following properties:
\smallskip \hskip 3.3mm {\bf (I)} The domain $D(A^*)$ is dense in $X^*$.

\smallskip \hskip 1.8mm {\bf (II)} The operator $A$ has a purely point
spectrum $\sigma $ which is either finite or has no finite limit points and
each $\lambda \in \sigma $ is of a finite multiplicity.

It is known \cite{Slmn87} , \cite{Weiss89}, etc. :

\begin{enumerate}
\item  for each $t\geq 0$ the operator $S\left( t\right) $ has a continuous
extension ${\cal S}\left( t\right) $on the space $V$ and the family of
operators ${\cal S}\left( t\right) :V\rightarrow V$ is the semigroup in the
class $C_{0}$ with respect to the norm of $V$ and the corresponding
infinitesimal generator ${\cal A}$ of the semigroup ${\cal S}\left( t\right) 
$ is the closed dense extension of the operator $A$ on the space $V$ with
domain $D\left( {\cal A}\right) $ $=X$;

\item  the set of the generalized eigenvectors of operators ${\cal A}$ $,$ $%
{\cal A}^{*}$ and $A,\,A^{*}\,$are the same$.$
\end{enumerate}

\smallskip \hskip 0.5mm {\bf (III)} There exists a time moment $T\ge 0$ such
that for all $v\in V$ and $t>T$ the function $x(t)={\cal S}(t{\cal )}v$ is
expanded in a series of generalized eigenvectors of the operator $A$,
converging with respect to the norm of $V$ for a certain grouping of terms
uniformly with respect to $t$ on an arbitrary interval $[T_{1},T_{2}]\ \
(T_{1}>T)$.

\section{Main results.}

Our main task consists of establishing conditions for independence of $K(t)$
at least for sufficiently large $t$.

\begin{definition}
A sequence $\{x_{i}\}_{i\in \hbox{\bb N}}$ of functions from $L_{2}^{{\rm loc%
}}[0,+\infty )$ is called {\it minimal} on $[0,\nu ]\ \ (\nu >0)$ if there
is a sequence $\{y_{j}\}_{j\in \hbox{\bb N}}$ of functions from $L_{2}[0,\nu
]$ such that 
\[
\int\limits_{0}^{\nu }\left( x_{i}(t)y_{j}(t)\right) \,dt=\delta _{ij}\qquad
(i,j\in \hbox{\bb N}) 
\]
where $\delta _{ij}$ is the Kronecker symbol. The sequence $\{y_{j}\}_{j\in %
\hbox{\bb N}}$ is called a sequence {\it biorthogonal} to the sequence $%
\{x_{j}\}_{j\in \hbox{\bb N}}\,$on $\left[ 0,\nu \right] $.
\end{definition}

Let the numbers $\lambda _{j}\in \sigma \ \ (j\in \hbox{\bb N})$ be
enumerated in the order of non-decreasing absolute values, let $\alpha _{j}$
be the multiplicity of $\lambda _{j}\in \sigma $, and let 
\[
\varphi _{jkl}\quad {\rm and\quad }\psi _{jkl,}\qquad j\in \hbox{\bb N};\
k=1,\ldots ,m_{j};\ l=1,2,\ldots ,\beta _{jk};\ \sum_{k=1}^{m_{j}}\beta
_{jk}=\alpha _{j}\ 
\]
be the generalized eigenvectors of the operators $A$ and $A^{*}$,
respectively, such that 
\begin{equation}
(\varphi _{jp\beta _{p}-l+1,}\psi _{ksq})=\delta _{jk}\delta _{ps}\delta
_{lq}  \label{e4}
\end{equation}
\[
j,k\in \hbox{\bb N};\ p=1,\ldots ,m_{j};\ l=1,\ldots ,\beta _{jp};\
s=1,\ldots ,m_{k};\ q=1,\ldots ,\beta _{ks}. 
\]

\begin{theorem}
\label{T1} If $t_{1}\le t_{2}$, then $K(t_{1})\subseteq K(t_{2})$. If the
properties (I)--(III) hold and the sequence the sequence $\{f_{jk}\}_{jk}$
of functions 
\begin{equation}
f_{jk}(t)=(-t)^{k}\exp (-\lambda _{j}t),\qquad j\in \hbox{\bb N};\
k=1,\ldots ,\alpha _{j};\ t\in [0,+\infty )  \label{e5}
\end{equation}
is minimal on $[0,\nu ]$, then

$\overline{K(t_{1})}=\overline{K(t_{2})}$ if $t_{1},t_{2}>T+\nu ${.}
\end{theorem}

{\bf Proof}\footnote{%
Some details of the proof from [Shklyar] are applicable for the case of
unbounded operator $B$ and can be omitted, however we repeat them here for
the sake of reader's convenience.}. A weak solution $x\left( t\right) $ of
equation (\ref{e1}) with the initial condition (\ref{e2}) is defined by the
following representation formula\cite{Lssk83},\cite{Slmn87}, \cite{Weiss89}

\begin{equation}
x\left( t,x_0,u\left( \cdot \right) \right) =S\left( t\right)
x_0+\int\limits_0^tS\left( t-\tau \right) Bu\left( \tau \right) d\tau .
\label{e6}
\end{equation}
Hence the attainable set $K\left( t\right) $ is defined by the formula

\begin{equation}
K_{t}=\left\{ x\in X:\exists u\left( \cdot \right) \in L_{2}\left( \left[
0,t_{1}\right] ,U\right) ,\,x=\int_{0}^{t}S\left( t-\tau \right) Bu\left(
\tau \right) d\tau \right\} ,  \label{e7}
\end{equation}
so one can prove the inclusion $K(t_{1})\subseteq K(t_{2})$ as well as in 
\cite{Shk95}

Let $P_j$ be a projector on the {\it generalized eugenspace} of $A$ at $%
\lambda _j\in \sigma \ \ (j\in \hbox{\bb N}),$ and let

\[
\Lambda _j=\left\{ 
\begin{array}{cccc}
{\lambda _j} & {1} & {\ldots } & {0} \\ 
{0} & {\lambda _j} & {\ldots } & {0} \\ 
0 & 0 & {\ldots } & 1 \\ 
0 & 0 & ... & {\lambda _j}
\end{array}
\right\} 
\]
be the Jordan $(\beta _j\times \beta _j)$-matrix. We have 
\begin{equation}
\left( S(t)P_jx,g\right) =(\Phi _j,g)\exp \left( \Lambda _jt\right) \left(
x,\Psi _j\right) ^{{\rm T}},\forall x\in X.\,,
\end{equation}
where 
\begin{eqnarray*}
\Phi _j &=&\left\{ \varphi _{jk1},\varphi _{jk2},\ldots ,\varphi _{jk\beta
_{jk}}\right\} \qquad \,\,\,\,\,\,\,k=1,\ldots ,m_j, \\
\Psi _j &=&\left\{ \psi _{jk1},\psi _{jk2},\ldots ,\psi _{jk\beta
_{jk}}\right\} \qquad \,\,\,\,\,\ k=1,\ldots ,m_j, \\
(\Phi _j,g) &=&\left\{ (\varphi _{jk1},g),\ldots ,(\varphi _{jk\beta
_{jk}},g)\right\} \qquad {k=1,\ldots ,}m{_j}, \\
(x,\Psi _j) &=&\left\{ (x,\psi _{jk1}),\ldots ,(x,\psi _{jk\beta
_{jk}})\right\} \qquad {k=1,\ldots ,}m.
\end{eqnarray*}

Now we will prove inclusion $K(t_{2})\subseteq K(t_{1})$ for all $T\le
t_{1}<t_{2}$. Let $g\in K(t_{1})^{\perp }$. By (\ref{e7}) 
\begin{equation}
\left( \int\limits_{0}^{t_{1}}S(t_{1}-\tau )Bu\left( \tau \right) d\tau
,g\right) \equiv 0,\,\forall u\in L\left( _{2}\left[ 0,t_{1}\right]
,U\right) .  \label{e7.2}
\end{equation}
If the linear operator $B$ is bounded, then

\begin{equation}
\left( \int\limits_{0}^{t_{1}}S\left( t_{1}-\tau \right) Bu\left( \tau
\right) d\tau ,g\right) =\int\limits_{0}^{t_{1}}\left( S\left( t_{1}-\tau
\right) Bu\left( \tau \right) ,g\right) d\tau .  \label{e8}
\end{equation}
If the linear operator is unbounded, then we cannot use (\ref{e8}) to
continue the proof as in [Shklyar], because there exists $u\left( \cdot
\right) \in L_{2}\left( \left[ 0,t_{1}\right] ,U\right) $ and $\tau \in
\left[ 0,t_{1}\right] $, such that $S\left( t-\tau \right) Bu\left( \tau
\right) \notin X$, and besides $S\left( t-\tau \right) Bu\left( \tau \right)
\in V$, but we cannot assure $g\in V^{*}$.

Let 
\[
u_1\left( t\right) =\left\{ 
\begin{array}{cc}
0, & t_1-T<t\leq t_1; \\ 
u\left( t\right) , & 0\leq t\leq t_1-T,
\end{array}
\right. 
\]
where $u\left( \cdot \right) \in L_2\left( \left[ 0,t_1-T\right] ,U\right) $%
. If follows from (\ref{e8}) that 
\begin{equation}
\left( \int\limits_0^{t_1-T}S(t_1-\tau )Bu\left( \tau \right) d\tau
,g\right) \equiv 0,\,\forall u\left( \cdot \right) \in L_2\left( \left[
0,t_1-T\right] ,U\right) .  \label{e9}
\end{equation}

Now let 
\[
u_2\left( t\right) =\left\{ 
\begin{array}{cc}
u\left( t\right) , & t_1-T<t\leq t_1; \\ 
0, & 0\leq t\leq t_1-T,
\end{array}
\right. 
\]
where $u\left( \cdot \right) \in L_2\left( \left[ t_1-T,t_1\right] ,U\right) 
$. Again, it follows from (\ref{e8}) that 
\begin{equation}
\left( \int\limits_{t_1-T}^{t_1}S(t_1-\tau )Bu\left( \tau \right) d\tau
,g\right) \equiv 0,\,\forall u\left( \cdot \right) \in L_2\left( \left[
t_1-T,t_1\right] ,U\right) .  \label{e9.1}
\end{equation}

The sets of the generalized eigenvectors of operators ${\cal A}$ $,$ ${\cal A%
}^{*}$ and $A,\,A^{*}\,$are the same, so in accordance with property (III)
we have

\begin{equation}
{\cal S}\left( t\right) v=\sum_{j=1}^\infty \Phi _j\exp \left( \Lambda
_jt\right) \left( v,\Psi _j\right) ^{\top },\forall t>T,\forall v\in V,
\label{e9.2}
\end{equation}
where $\sum_{j=1}^\infty $ is considered with respect to the norm of $V$.

Consider the sequence $S_n\left( t\right) v$ of partial sums of series (\ref
{e9.2}) 
\begin{equation}
{\cal S}_n\left( t\right) v=\sum_{j=1}^n\Phi _j\exp \left( \Lambda
_jt\right) \left( v,\Psi _j\right) ^{{\rm T}}.  \label{e9.3}
\end{equation}
Obviously $S_n\left( t\right) v\in X,\forall v\in V$.

One can show that 
\begin{equation}
\left( \Phi _{j},g\right) \exp \left( \Lambda _{j}t\right) \left( v,\Psi
_{j}^{{\rm T}}\right) =\sum_{k=0}^{\beta _{j}}\exp \left( \lambda
_{j}t\right) \left( \Phi _{j},g\right) \frac{t^{k}}{k!}E_{j}^{k}\left(
v,\Psi _{j}^{{\rm T}}\right) ,\   \label{e9.3.1}
\end{equation}
where $\beta _{j}\times \beta _{j}$-matrix $E_{j}$ is defined by $\,$%
\[
E_{j}=\left\{ 
\begin{array}{cccc}
{0} & {1} & {\ldots } & {0} \\ 
{0} & {0} & {\ldots } & {0} \\ 
0 & 0 & {\ldots } & 1 \\ 
0 & 0 & ... & {0}
\end{array}
\right\} . 
\]

Denote by $g_n\left( v,t\right) $ the linear functional (with respect to $v$)

\begin{equation}
g_{n}\left( v,t\right) =\left( {\cal S}_{n}\left( t\right) v,g\right)
=\sum_{j=1}^{n}\left( \Phi _{j},g\right) \exp \left( \Lambda _{j}t\right)
\left( v,\Psi _{j}^{{\rm T}}\right) .  \label{e9.4}
\end{equation}

It follows from above considerations that for any $j=1,2,...$ the linear
operator acting from $V\,$\thinspace to \hbox{\bb R}$^{\beta _{j}}\,$
defined by $\left( v,\Psi _{j}^{{\rm T}}\right) $ is bounded. Hence for any
natural $n$ the functional $g_{n}\left( v,t\right) \,\,$ is a linear bounded
functional, and it follows from (\ref{e9}) and (\ref{e9.2}), that for $%
v=Bu\left( t_{1}-\tau \right) $, where $u\left( \tau \right) =0$ for $\tau
\in (t_{1}-T,t_{1}]\,$

\begin{equation}
\lim_{n\rightarrow \infty }\int\limits_{0}^{t_{1}-T}g_{n}\left( Bu\left(
t_1-\tau \right) ,t_1-\tau \right) d\tau =0.  \label{e9.5}
\end{equation}

Let $w\in U\,\,$and $\gamma _{kl}\left( t\right) ,k=1,2,...,l=1,2,...,\beta
_{k}\,\,$be the sequence of functions biorthogonal to the sequence (\ref{e6}%
) on $\left[ 0,t_{1}-T\right] .$ $\,$Substituting (\ref{e9.4}) to (\ref{e9.5}%
) and using $u(t)=w\gamma _{kl}\left( t\right) $ we obtain after
computations 
\begin{equation}
\left( \Phi _{j},g\right) E_{j}^{k}B^{*}\Psi _{j}^{{\rm T}%
}=0,j=1,2,...,\,k=1,...,\beta _{j},  \label{e9.6}
\end{equation}
On account of (\ref{e9.6}), (\ref{e9.4}) and (\ref{e9.3.1}) we obtain 
\[
\lim_{n\rightarrow \infty }\int\limits_{0}^{t_{2}-T}g_{n}\left( Bu\left(
t_{2}-\tau \right) ,t_{2}-\tau \right) d\tau =0, 
\]
so 
\begin{equation}
\left( \int\limits_{0}^{t_{2}-T}S(t_{2}-\tau )Bu\left( \tau \right) d\tau
,g\right) \equiv 0,\,\forall u\left( \cdot \right) \in L_{2}\left( \left[
0,t_{2}-T\right] ,U\right) .  \label{e9.7}
\end{equation}
Joining (\ref{e9.7}) and (\ref{e9.1}), we obtain 
\begin{equation}
\left( \int\limits_{0}^{t_{2}}S(t_{2}-\tau )Bu\left( \tau \right) d\tau
,g\right) \equiv 0,\,\forall u\in L_{2}\left( \left[ 0,t_{1}\right]
,U\right) .  \label{e10}
\end{equation}

The latter identity imply the inclusion $g\in K(t_2)^{\perp }$. Thus, $%
K(t_1)^{\perp }\subseteq K(t_2)^{\perp }$. Hence $\overline{K(t_2})\subseteq 
\overline{K(t_1})$. Since $K(t_1)\subseteq K(t_2)$ for all $t_1$ with $%
t_1<t_2$, we obtain $\overline{K(t_1})=\overline{K(t_2})$ for all $t_1$ and $%
t_2$ with $T+\nu <t_1<t_2$. This proves the theorem.

\subsection{Controllability conditions.}

Theorem \ref{T1} can be applied for various control problems.

In this section we will show how the proof of Theorem \ref{T1} provides a
possibility to obtain an approximate null-controllability criterion for the
abstract control problem with unbounded input operator. We will consider
this kind of controllability only, but other kinds of controllability can be
investigated also.

\subsubsection{Approximate null-controllability conditions for equation
(\ref{e1}.}

\smallskip 
Denote 
\[
{\rm Range}\{\lambda I-A,R_{\mu }B\}=\{z\in X:\exists x\in X,\,\exists u\in
U,\,{\rm \,}z=(\lambda I-A)x+R_{\mu }Bu\} 
\]

\begin{theorem}
\label{T3} Let $\mu \notin \sigma .$ If the properties (I)--(III) hold and
the sequence the sequence (\ref{e5}) is minimal on $[0,\nu ]$, then $\,$for
equation (\ref{e1}) to be approximately null-controllable on $[0,t_{1}]$ it
is necessary and, for $t_{1}>T+\nu $, sufficient that 
\begin{equation}
\overline{{\rm Range}\{\lambda I-A,R_{\mu }B\}}=X,\,\forall \lambda \in
\sigma  \label{e13}
\end{equation}
\end{theorem}

{\bf Proof. } {\it Sufficiency}. We obtained above that $g\in
K(t_{1})^{\perp }$ provided $t_{1}>T$ implies the identity (\ref{e9.6}). One
can easy see that the condition (\ref{e13}) is equivalent to the condition 
\begin{equation}
B^{*}R_{\mu }^{*}\psi _{j\beta _{j}}\neq 0,\,j\in \hbox{\bb N}  \label{e13.1}
\end{equation}
Here $\psi _{j\beta _{j}}\,$ is the eigenvector of the adjoint operator $%
A^{*}$ corresponding to the eigenvalue $\lambda _{j}\in \sigma $. Since the
eigenvectors of the operators $A^{*}\,$and ${\cal A}^{*}\,$are the same , we
have $\psi _{j\beta _{j}}\in V^{*},$ so $B^{*}\psi _{j\beta _{j}}\,$is
well-defined and it follows from (\ref{e13.1}) that 
\begin{equation}
B^{*}R_{\mu }^{*}\psi _{j\beta _{j}}=\left( \mu -\lambda \right) B^{*}\psi
_{j\beta _{j}}\neq 0,\,j\in \hbox{\bb N},\,\mu \notin \sigma ,\,\lambda \in
\sigma .  \label{e13.2}
\end{equation}
Hence (\ref{e13.2}) yields 
\begin{equation}
B^{*}\psi _{j\beta _{j}}\neq 0,\,j\in \hbox{\bb N}.  \label{e13.3}
\end{equation}

Solving the linear algebraic system (\ref{e9.6}) provided that (\ref{e13.3})
holds, we obtain 
\begin{equation}
(\Phi _{j},g)=0\qquad (j\in \hbox{\bb N}).  \label{e14}
\end{equation}
This and property (III) imply $S^{*}(t_{1})g=0$, therefore, $g\in \,{\rm %
Range}S(t_{1})^{\perp }$. We have $K(t_{1})^{\perp }\subseteq {\rm Range}%
\,S(t_{1})^{\perp }$, hence $\overline{{\rm Range}\,S(t_{1})}\subseteq 
\overline{K(t_{1}})$. The latter relation is equivalent to the approximate
null-controllability of equation (1) on $[0,t_{1}]$. This proves the
sufficiency of (\ref{e13}).

\smallskip {\it Necessity}. If condition (\ref{e13}) does not hold, then
there exists $\lambda \in \sigma $ and $g\in X^{*},g\neq 0$ such that 
\begin{eqnarray}
(\lambda x-Ax,g) &=&0,\quad \forall x\in D(A),  \label{e14.01} \\
\left( R_{\mu }Bu,g\right) &=&0,\,\,\,\,\,\forall u\in U.  \label{e14.011}
\end{eqnarray}
It follows from (\ref{e14.01}) that the vector $g$ is the eigenvector of the
operator $A^{*}$ corresponding to eigenvalue$\,\lambda $. Since the
eigenvectors of the operators $A^{*}\,$and ${\cal A}^{*}\,$are the same we
have $g\in V^{*}$, the scalar product $\left( Bu,g\right) \,$is
well-defined. and one can write (\ref{e14.011}) in the form 
\begin{equation}
\left( R_{\mu }Bu,g\right) =\left( Bu,R_{\mu }^{*}g\right) =\left( Bu,\left( 
\bar{\mu}-\bar{\lambda}\right) g\right) =\left( \mu -\lambda \right) \left(
Bu,g\right) .  \label{e14.012}
\end{equation}
It follows from (\ref{e14.01})--(\ref{e14.012}) that$\left( S(t)Bu,g\right)
\equiv 0$ for all $t\in [0,+\infty )$ and $u\in U$, but $S^{*}(t)g\neq 0$
for all $t\in [0,+\infty )$. Hence, $g\in K(t_{1})^{\perp }$, but $g\notin 
{\rm Range}\,S(t_{1})^{\perp }$. This proves the necessity of (\ref{e13}).

If the operator $A$ is not self adjoint, then it is not trivial problem to
calculate the adjoint operator $A.$ If the operator $A^{*}$ is calculated
then one can use instead condition (\ref{e13}) one of the conditions:

\begin{enumerate}
\item  {\it for any }$\lambda \in \sigma $ and any $\mu \notin \sigma $ {\it %
the system of equations with respect to } $g\in X^{*}$%
\begin{eqnarray}
\lambda g-A^{*}g=0,  \label{e14.02} \\
B^{*}R_{\mu }^{*}g=0  \nonumber
\end{eqnarray}
{\it has only trivial solution;}

\item  {\it for any }$\lambda \in \sigma $ {\it the system of equations with
respect to } $g\in V^{*}$%
\begin{eqnarray}
\lambda g-A^{*}g=0,  \label{e14.021} \\
B^{*}g=0  \label{14.022}
\end{eqnarray}
{\it has only trivial solution;}
\end{enumerate}

By (\ref{e14.01})-(\ref{e14.021}) one can easy obtain rank conditions for
approximate null-controllability of equation (\ref{e1}).

\begin{remark}
Since the generalized eigenvectors of the operators $A^{*}\,$and ${\cal A}%
^{*}\,$are the same we have $B^{*}\Psi _{j}$ to be well-defined,\thinspace
\thinspace $j=1,2,...$ .
\end{remark}

\begin{theorem}
\label{T6}Let $I$ be $\beta _{j}\times \beta _{j}$ unit $\,$matrix. If the
properties (I)--(III) hold and the sequence the sequence (\ref{e5}) is
minimal on $[0,\nu ]$, then $\,$for equation (\ref{e1}) to be approximately
null-controllable on $[0,t_{1}]$ it is necessary and, for $t_{1}>T+\nu $,
sufficient that 
\begin{equation}
{\rm rank}\left\{ \lambda _{j}I_{j}-\Lambda _{j},B^{*}\Psi _{j}\right\}
=\beta _{j},\,j=1,2,...\,.  \label{e14.03}
\end{equation}
\end{theorem}

{\bf Proof. }It follows from (\ref{e14.021}) that there exists a vector $%
\eta \in {\hbox{\bb R}^{\beta _{j}}}$ such that $g=\Psi _{j}\eta $. Since $%
A^{*}\Psi _{j}=\Psi _{j}\Lambda _{j}^{*},\,j=1,2,...,$ we obtain from (\ref
{e14.021}), that for any $j=1,2,...$%
\begin{eqnarray}
\Psi _{j}\left( \lambda _{j}I_{j}-\Lambda _{j}^{*}\right) \eta &=&0,
\label{e14.04} \\
B^{*}\Psi _{j}\eta &=&0.  \label{e14.05}
\end{eqnarray}

The condition 
\[
{\rm rank}\Psi _{j}=\beta _{j},\,j=1,2,... 
\]
yields the equivalence between (\ref{e14.04})-(\ref{e14.05}) and 
\begin{eqnarray}
\left( \lambda _{j}I_{j}-\Lambda _{j}^{*}\right) \eta &=&0,  \label{e14.06}
\\
B^{*}\Psi _{j}\eta &=&0.  \label{e14.07}
\end{eqnarray}
Thus (\ref{e14.021}) is equivalent to (\ref{e14.06})-(\ref{e14.07}), where $%
g=\Psi _{j}\eta $,\thinspace so (\ref{e14.06})-(\ref{e14.07}) holds if and
only if $\eta =0$. This shows the validity of (\ref{e14.03}).

\subsubsection{Approximate null-controllability conditions for Abstract
Boundary Control Problem.}

Let $X,$ $U,$ $Y$ be Banach spaces. Consider the abstract boundary control
problem

\begin{eqnarray}
\dot{x}\left( t\right) &=&Lx\left( t\right) ,  \label{e14.1} \\
Gx\left( t\right) &=&Bu\left( t\right) ,  \label{e14.2} \\
x\left( 0\right) &=&x_{0},  \label{e14.3}
\end{eqnarray}
where $L:X\rightarrow X$ is a linear unbounded operator with dense domain $%
Z=D\left( L\right) ,$ $B:U\rightarrow Y$ is a linear bounded injective
operator, $G:Z\rightarrow Y$ is a linear bounded operator satisfying the
following conditions:

\begin{itemize}
\item  $G$ is onto, {\rm Ker}$G$ is dense in $X$ ;

\item  there exists a $\mu \in {\hbox{\bb R}}$ such that $\mu I-L$ is onto
and {\rm Ker}$\left( \mu I-L\right) \bigcap $ {\rm Ker}$G=\emptyset $.
\end{itemize}

The problem (\ref{e14.1})-(\ref{e14.3}) is assumed to be well-posed. This
problem is the abstract model for classical control problem described by
linear partial differential equations of parabolic type when the control
acts through the boundary and if the measurement can only be realized at a
few points of the corresponding spatial domain. The conditions (\ref{e14.2})
can be considered as abstract boundary conditions.

Now we transform the abstract boundary control problem (\ref{e14.1})-(\ref
{e14.3}) to the problem (\ref{e1})-(\ref{e2}). Consider the space $W=$ {\rm %
Ker}$G.$ We have $W\subset Z\subset X$ with continuous dense injection. The
operator $A:W\rightarrow X$ is defined by 
\begin{equation}
Ax=Lx\,\,{\rm for\,\,}x\in W.  \label{e14.4}
\end{equation}
For any $y\in Y$ define 
\begin{equation}
\hat{B}y=Lx-Ax,\,x\in G^{-1}\left( y\right) =\{z\in Z:Gx=y\}.  \label{e14.5}
\end{equation}
The operator $\hat{B}:Y\rightarrow V$ defined by (\ref{e14.5}) is a bounded
operator, but it is unbounded as an operator $\hat{B}:Y\rightarrow X.$ Given 
$u\in U$ denote $\tilde{B}u=\hat{B}Bu.$ The operator $\tilde{B}:U\rightarrow
V$ is bounded, but the corresponding operator $\tilde{B}:U\rightarrow X$ is
unbounded. It follows from (\ref{e14.5}) that 
\begin{eqnarray}
Lx &=&Ax+\tilde{B}u,  \label{e14.6} \\
Gx &=&Bu.
\end{eqnarray}

Since the abstract boundary control problem under consideration is uniformly
well-posed the operator $A$ generates the strongly continuous semigroup of
bounded operators in the class $C_0.$ Hence given abstract boundary control
problem is equivalent to the control problem 
\begin{eqnarray}
\dot x\left( t\right) &=&Ax\left( t\right) +\tilde Bu\left( t\right) ,
\label{e14.7} \\
x\left( 0\right) &=&x_0.  \label{e.14.8}
\end{eqnarray}
So it follows from above considerations and Theorem that

\begin{theorem}
\label{T4} Let $\mu \notin \sigma .$ If the properties (I)--(III) hold and
the sequence (\ref{e5}) is minimal on $[0,\nu ]$, then for equation (\ref
{e14.1})-(\ref{e14.2}) to be approximately null-controllable on $[0,t_{1}]$
it is necessary and, for $t_{1}>T+\nu $, sufficient that 
\begin{equation}
\overline{{\rm Range}\{\lambda I-A,R_{\mu }\hat{B}B\}}=X,\,\forall \lambda
\in \sigma .  \label{e14.9}
\end{equation}
\end{theorem}

Together with equation (\ref{e14.1})-(\ref{e14.3}) consider the abstract
elliptic equation 
\begin{equation}
Lx=\mu x  \label{e14.10}
\end{equation}
\begin{equation}
Gx=y  \label{e14.11}
\end{equation}

Since the problem (\ref{e14.1})-(\ref{e14.3}) is uniformly well-posed then
for any $y\in Y$ there exists the solution $x_\mu =D_\mu y$ of the equation (%
\ref{e14.10})-(\ref{e14.11}), where $D_\mu :Y\rightarrow X$ is a linear
bounded operator. It follows from (\ref{e14.10})-(\ref{e14.11}) and (\ref
{e14.6}) that 
\[
Ax+\hat By=\mu x. 
\]
Hence 
\[
x=R_\mu \hat By 
\]
so 
\begin{equation}
D_\mu =R_\mu \hat B.  \label{e14.12}
\end{equation}
Using (\ref{e14.12}) in (\ref{e14.9}) we obtain that the following theorem
is valid:

\begin{theorem}
\label{T5}Let $\mu \notin \sigma .$ If the properties (I)--(III) hold and
the sequence (\ref{e5}) is minimal on $[0,\nu ]$, then for equation (\ref
{e14.1})-(\ref{e14.2}) to be approximately null-controllable on $[0,t_{1}]$
it is necessary and, for $t_{1}>T+\nu $, sufficient that 
\begin{equation}
\overline{{\rm Range}\{\lambda I-A,D_{\mu }B\}}=X,\,\,\forall \lambda \in
\sigma .  \label{e14.13}
\end{equation}
\end{theorem}

There are a lot of methods to calculate the operator $D_{\mu }$ for a given
concrete boundary control problem [Butkovskii].

\section{Examples.}

In this section we will apply the results obtained in previous sections to a
general class of linear neutral functional differential equations, which don$%
^{\prime }$t fit into the framework of \cite{Shk95}.

\subsection{General neutral functional differential equations.}

Consider the linear neutral functional differential equation\cite{Hale77} 
\begin{eqnarray}
&&\ \ \frac{d}{dt}\left( x\left( t\right) -\int\limits_{-h}^{0}dA_{0}(\tau
)x(t+\tau )-\int\limits_{-h}^{0}dB_{0}(\tau )u(t+\tau )\right)  \label{e15}
\\
&=&\ \int\limits_{-h}^{0}dA(\tau )x(t+\tau )+\int\limits_{-h}^{0}dB(\tau
)u(t+\tau ),
\end{eqnarray}
with initial conditions 
\begin{equation}
x\left( 0\right) -\int\limits_{-h}^{0}dA_{0}(\tau )\varphi _{1}(\tau
)-\int\limits_{-h}^{0}dB_{0}(\tau )\varphi _{2}(\tau )=\varphi ^{0},x\left(
\tau \right) =\varphi _{1}(\tau ),u(\tau )=\varphi _{2}(\tau ),
\label{e15.1}
\end{equation}
where $x\left( t\right) \in {\hbox{\bb R}}^{n},\,\,\varphi _{1}(\cdot )\,\in
L_{2}\left( \left[ -h,0\right] ,{\hbox{\bb R}}^{n}\right) ,\,\varphi
_{2}(\cdot )\,\in L_{2}\left( \left[ -h,0\right] ,{\hbox{\bb R}}^{n}\right)
,\,u\left( t\right) \in {\hbox{\bb R}}^{r};\,\,A_{0}(\cdot ),\,B_{0}(\cdot
),\,A(\cdot ),$ $B(\cdot )$ are $(n\times n)$ and $(n\times r)$%
-matrix-functions of bounded variation, respectively, and the
matrix-function $A_{0}\left( \tau \right) $ satisfies the condition\footnote{%
The condition (\ref{e15.2}) provides the existence and uniqueness for
solutions of (\ref{e15})-(\ref{e15.1}) [Heil].} 
\begin{equation}
A_{0}\left( 0\right) =\lim_{\tau \rightarrow 0^{+}}A_{0}\left( \tau \right) .
\label{e15.2}
\end{equation}
The state space of the equation under consideration is 
\[
X={\hbox{\bb R}}^{n}\times L_{2}\left( [-h,0],{\hbox{\bb R}}^{n}\right)
\times L_{2}\left( [-h,0],{\hbox{\bb R}}^{r}\right) \footnote{%
The state spaces of [Salamon] is the Hilbert space 
\[
{\hbox{\bb R}}^n\times L_2\left( \left[ -h,0\right] ,{\hbox{\bb R}}^n\right)
\times L_2\left( \left[ -h,0\right] ,{\hbox{\bb R}}^r.\right) 
\]
}. 
\]

Let $x_{t}$ be the function defined by 
\[
x_{t}\left( \tau \right) =x\left( t+\tau \right) ,\,u_{t}\left( \cdot
\right) =u\left( t+\tau \right) ,-h\leq \tau \leq 0. 
\]

Following \cite{Slmn87}, we will describe equation (\ref{e15})-(\ref{e15.1})
by well-posed abstract boundary control problem 
\begin{eqnarray}
\frac{d}{dt}\bar{x}\left( t\right) &=&L\bar{x}\left( t\right) ,  \label{e16}
\\
G\bar{x}\left( t\right) &=&u\left( t\right) ,  \label{e16.1}
\end{eqnarray}
where$\,\,\,\,\,\,\,\,{Z}{}{=}$ 
\begin{equation}
{}{}\left\{ 
\begin{array}{c}
\bar{x}\in X:\varphi ^{0}\in {\hbox{\bb R}}^{n},\varphi _{1}\in
W^{1,2}\left( [-h,0],{\hbox{\bb R}}^{n}\right) ,\varphi _{2}\in
W^{1,2}\left( [-h,0],{\hbox{\bb R}}^{r}\right) \\ 
\varphi ^{0}=\varphi _{1}\left( 0\right) -\int\limits_{-h}^{0}dA_{0}(\tau
)\varphi _{1}(\tau )-\int\limits_{-h}^{0}dB_{0}(\tau )\varphi _{2}(\tau ),
\end{array}
\right\}  \label{e17}
\end{equation}

$\bar x\left( t\right) =\left( x\left( t\right) -\int\limits_{-h}^0dA_0(\tau
)x(t+\tau )-\int\limits_{-h}^0dB_0(\tau )u(t+\tau ),x_t,u_t\right) ,$

$\bar x\left( 0\right) =\bar x=\left( \varphi ^0-\int\limits_{-h}^0dA_0(\tau
)\varphi _1(\tau )-\int\limits_{-h}^0dB_0(\tau )\varphi _2(\tau ),\varphi
_1\left( \cdot \right) ,\varphi _2\left( \cdot \right) \right) ;$%
\begin{eqnarray}
L\bar x &=&\left( \int\limits_{-h}^0dA(\tau )\varphi _1(\tau
)-\int\limits_{-h}^0dB(\tau )\varphi _2(\tau ),\dot \varphi _1(\cdot ),\dot
\varphi _2(\cdot )\right) ,  \label{e18} \\
G\bar x &=&\dot \varphi _2(0)  \nonumber
\end{eqnarray}
for which the corresponding operator $A$ satisfies the conditions (I) - (II)
[10], and condition (III) holds for a wide class of (\ref{e15}) (see [14: p.
101]). For example, condition (III) holds for 
\[
A(\tau )=\sum_{j=0}^mA_j\chi _{[-h_{j+1},-h_j]}(\tau )\qquad \big(%
0=h_0<h_1<\ldots <h_m=h,\ -h\le \tau \le 0\big)
\]
\[
A_0(\tau )=\sum_{j=0}^mA_{0j}\chi _{[-h_{j+1},-h_j]}(\tau )\qquad \big(%
0=h_0<h_1<\ldots <h_m=h,\ -h\le \tau \le 0\big)
\]
where $A_j$ are $(n\times n)$-matrices and $\chi _{[-h_{j+1},-h_j]}(\tau )$
is the characteristic function of the interval $[-h_{j+1},-h_j).$

Let 
\[
\Delta (z)=\det \left\{ zI-\int\limits_{-h}^{0}dA_{0}(\tau )z\exp {z\tau }%
-\int\limits_{-h}^{0}dA(\tau )\exp {z\tau }\right\} . 
\]
Denote by $\omega $ the exponential of the function $\Delta (z)$ \cite{Lvn64}%
, i.e. 
\[
\omega =\overline{\lim_{|x|\rightarrow \infty }}{\frac{1}{|z|}}\log |\det {%
\Delta (z)}|. 
\]

\begin{lemma}
\label{L1} The sequence $(6)$ is minimal on $[0,\nu ]$ for any $\nu >\omega $%
.
\end{lemma}

{\bf Proof}. The assertion of the lemma has been proven in [Shklyar] for
hereditary case and without delays in control ($A_{0}\left( \tau \right) $
and $B_{0}\left( \tau \right) $ are constant on the segment $\left[
-h,0\right] ;\,B\left( \tau \right) =\left\{ 
\begin{tabular}{l}
$B,\,\,{\rm if\,\,}\tau =0,$ \\ 
$0,\,\,{\rm if\,\,}-h<\tau \leq 0.$%
\end{tabular}
\right. $ ).

One can use the same proof in the general case. There is only one difference:

\[
\Delta (z)=\sum_{j=0}^{n}r_{j}(z)z^{j}, 
\]
where $r_{j}(z)$ is represented as a finite sum of products of numbers 
\[
\int\limits_{-h}^{0}da_{jk}(\tau )\exp (-z\tau )\,\,{\rm and\,\,}%
\int\limits_{-h}^{0}da_{jk}^{0}(\tau )\exp (-z\tau )\qquad (j=1,\ldots ,n) 
\]
with $a_{jk}(\tau )\,$ and $a_{jk}^{0}(\tau )$ being the elements of the
matrix $A(\tau )$ and $A_{0}\left( \tau \right) \,$ correspondingly .
Further we finish the proof as well as in \cite{Shk95}.

\subsection{Partial hyperbolic systems.}

Consider the wave equation on $\left[ 0,t_{1}\right] \times \left[ 0,\pi
\right] $\cite{Weiss89}$:$%
\begin{eqnarray}
\frac{\partial ^{2}}{\partial t^{2}}w\left( t,\theta \right) &=&\frac{%
\partial ^{2}}{\partial \theta ^{2}}w\left( t,\theta \right) ,  \label{e19}
\\
w\left( t,0\right) &=&u\left( t\right) ,w\left( t,\pi \right) =u\left(
t\right) ,  \nonumber \\
w\left( 0,\theta \right) &=&\varphi _{0}\left( \theta \right) ,\frac{%
\partial }{\partial t}w\left( 0,\theta \right) =\varphi _{1}\left( \theta
\right) .  \nonumber
\end{eqnarray}

We assume the weak solution $w\left( t,\cdot \right) \in AC\left[ 0,\pi
\right] ,$where $AC\left[ 0,\pi \right] $\thinspace is\thinspace
the\thinspace space of absolutely continuous functions defined on $\left[
0,\pi \right] ;\,\frac{\partial }{\partial t}w\left( t,\cdot \right)
,\,\varphi _{0}\left( \cdot \right) $ and $\varphi _{1}\left( \cdot \right)
\in L_{2}\left[ 0,\pi \right] ;\,u\left( \cdot \right) \in L_{2}\left[
0,t_{1}\right] .$ Define $v\left( t,\cdot \right) =\frac{\partial }{\partial
t}w\left( t,\cdot \right) $%
\[
x=\left\{ 
\begin{tabular}{l}
$w\left( t,\cdot \right) $ \\ 
$v\left( t,\cdot \right) $%
\end{tabular}
\right\} 
\]

\[
X=W_{0}^{1}\left[ 0,\pi \right] \times L_{2}\left[ 0,\pi \right] , 
\]
where 
\begin{equation}
W_{0}^{1}\left[ 0,\pi \right] =\left\{ x\in AC\left[ 0,\pi \right] :x\left(
0\right) =x\left( \pi \right) =0,\frac{d}{d\theta }x\left( \cdot \right) \in
L_{2}\left[ 0,\pi \right] \right\} .  \label{e20}
\end{equation}
We have 
\begin{eqnarray}
\dot{x}\left( t\right) &=&Lx\left( t\right) ,  \label{e21} \\
Gx\left( t\right) &=&bu\left( t\right) ,  \nonumber
\end{eqnarray}
where the operators $A,G\,$ and $b$ are defined by 
\[
L=\left( 
\begin{array}{cc}
0 & 1 \\ 
\frac{\partial ^{2}}{\partial \theta ^{2}} & 0
\end{array}
\right) , 
\]
\[
Gx=\left( 
\begin{array}{l}
w\left( 0\right) \\ 
w\left( \pi \right)
\end{array}
\right) , 
\]
\[
b=\left( 
\begin{array}{l}
1 \\ 
1
\end{array}
\right) 
\]
with domain

\[
D\left( A\right) =\left( W_{0}^{1}\left[ 0,\pi \right] \bigcap W^{2}\left[
0,\pi \right] \right) \times W_{0}^{1}\left[ 0,\pi \right] ,D\left( G\right)
=C\left( \left[ 0,\pi \right] ,\hbox{\bb R}^{2}\right) ,
\]
where 
\begin{equation}
W^{2}\left[ 0,\pi \right] =\left\{ x\in AC\left[ 0,\pi \right] :\frac{d}{%
d\theta }x\left( \cdot \right) \in AC\left[ 0,\pi \right] ,\,\frac{d^{2}}{%
d\theta ^{2}}x\left( \cdot \right) \in L_{2}\left[ 0,\pi \right] \right\} .
\label{e22}
\end{equation}
The above operator $A$ generates a contraction group $S\left( t\right) $ on $%
X$.$\ $

The spectrum $\sigma $ of the operator $A\ $is defined by 
\[
\sigma =\left\{ \lambda \in \hbox{\bb C}\right\} :\lambda =\pm
ki,k=0,1,2,... 
\]
$\ $

We have the sequence of functions 
\[
e^{kt},k=0,\pm 1,\pm 2,... 
\]
be minimal on $\left[ 0,2\pi \right] ,$and all the properties ({\bf I)-(III) 
}for any $T>0$. hence the next theorem follows from Theorem \ref{T1}:

\begin{theorem}
\label{T7}\label{T8} The closure of the attainable set for equation (\ref
{e19}) doesn$^{\prime }$t depend on $t$ for any $t>2\pi $.
\end{theorem}

\section{Concluding remarks.}

Properties of the attainable set $K(t)$ for equation (\ref{e1}) with
unbounded input operators are considered. Results of \cite{Shk95} have been
generalized for such classes of abstract evolution equations. Based on
Theorem \ref{T1} the null-controllability criterion for equation (\ref{e1})
is obtained. Application to general functional differential equation of
neutral type and hyperbolic systems have been considered. As well as in \cite
{Shk95}, property (III) and the minimality of the functions (\ref{e6})
provide the required independence of $t$ for $\overline{K(t)}$.

By duality principle one can obtain observability conditions for abstract
evolution equation with unbounded output operators.



\begin{thebibliography}{9}
\bibitem{DPrt84}  Da Pratto, G.: {\it Abstract differential equations and
extrapolation spaces. }Lecture Notes in Mathematics, 1184, Springer-Berlag,
Berlin, New York, 1984.

\bibitem{Hale77}  Hale, J.: {\it Theory of Functional Differential Equations.%
} Berlin - Heidelberg - New York: Springer-Verlag 1977. \smallskip 

\bibitem{Lssk83}  Lasieska, A.: {\it Unified theory for abstract parabolic
boundary problems -- A semigroup approach. }Appl. Math. Optimiz. , 10
(1983), 225 - 286.

\bibitem{Lvn64}  Levin, B. : {\it Distribution of Zeroes of Entire Functions.%
} Providence (R.I.): Amer. Math. Soc. 1964. \smallskip 

\bibitem{Ngl86}  Nagel, R.: {\it One-parameter semigroups of positive
operators. }Lecture Notes in Notes in Mathematics, 1184, Springer-Berlag,
Berlin, New York, 1984.

\bibitem{Slmn87}  Salamon, D.: {\it Infinite dimensional linear systems with
unbounded control and observation: a functional analytic approach.} Trans.
Amer. Math. Soc. 300 (1987), 383 - 431. \smallskip 

\bibitem{Shk91}  Shklyar, B.: {\it Controllability of linear systems with
distributed parameters.} Diff. Equ. 26 (1991), 326 - 335. \smallskip

\bibitem{Shk95}  Shklyar, B.: {\it On attainable set for abstract control
problem with application to controllability}. Z.Angew Math. und Mech, 1995.

\bibitem{Weiss89}  Weiss, G.: {\it Admissibility of unbounded control
operators.. }SIAM J. Contr. and Optimiz. 27 (1989), 527 - 545.
\end{thebibliography}
\end{document}